\documentclass[12pt]{amsart}
\usepackage{amsfonts,color,latexsym,amssymb,amsmath,epsfig,rotating,array}

\newtheorem{dfn}{Definition}[section]

\newtheorem{thm}[dfn]{Theorem}

\newtheorem{conj}[dfn]{Conjecture}

\definecolor{purple}{rgb}{.5,0,.5}
\definecolor{red}{rgb}{.6,0,0} 
\definecolor{green}{rgb}{0,.5,0} 

\renewcommand{\qed}{$\blacksquare$}

\makeatletter
\renewcommand*\env@matrix[1][c]{\hskip -\arraycolsep
  \let\@ifnextchar\new@ifnextchar
  \array{*\c@MaxMatrixCols #1}}
\makeatother

\newlength{\burg}
\newlength{\koi}
\newlength{\sma}
\newlength{\jmr}
\newlength{\belpa}
\newlength{\taovu}
\newlength{\hwl}
\newlength{\khov}
\newcommand{\thth}{^{\text{\underline{th}}}} 

\newcommand{\rd}{^{\text{\underline{rd}}}}

\newcommand{\eps}{\varepsilon}
\newcommand{\F}{\mathbb{F}}
\newcommand{\Q}{\mathbb{Q}}

\newcommand{\C}{\mathbb{C}}

\newcommand{\N}{\mathbb{N}}
\newcommand{\Z}{\mathbb{Z}}

\newcommand{\cO}{\mathcal{O}}

\newcommand{\rmv}[1]{}

\author{Qi Cheng} 
\email{qcheng@cs.ou.edu}
\address{School of Computer Science, 
University of Oklahoma, Norman, OK \ 73019} 
\author{Shuhong Gao} 
\email{sgao@math.clemson.edu} 
\address{Department of Mathematical Sciences, 
Clemson University, Clemson, SC \ 29634-0975} 
\author{J.\ Maurice Rojas}
\email{rojas@math.tamu.edu} 
\address{TAMU 3368, College Station, TX \ 77843-3368} 
\thanks{Partially supported by NSF grant CCF-1409020 
and the American Institute of Mathematics. This work 
was also supported by LABEX MILYON (ANR-10-LABX-0070) of Universit\'e de 
Lyon, within the program ``Investissements d'Avenir'' (ANR-11-IDEX-0007) 
operated by the French National Research Agency (ANR)}  
\author{Daqing Wan} 
\email{dwan@math.uci.edu} 
\address{Department of Mathematics, University of California, 
Irvine, CA\  92697-3875 }
\title[Sparse Univariate Polynomials Over Finite Fields] 
{\mbox{}\\
\vspace{-1in} 
Sparse Univariate Polynomials \\ with Many Roots Over Finite Fields} 

\begin{document}

\begin{abstract} 
Suppose $q$ is a prime power and $f\!\in\!\F_q[x]$ is a univariate polynomial 
with exactly $t$ monomial terms and degree $<\!q-1$.
To establish a finite field analogue of Descartes' Rule, Bi, Cheng, and Rojas 
(2013) proved an upper bound of $2(q-1)^{\frac{t-2}{t-1}}$ on the number of 
cosets in $\F^*_q$ needed to cover the roots of $f$ in $\F^*_q$. Here, we  
give explicit $f$ with root structure approaching this bound: 
For $q$ a $(t-1)$-st power of a prime we give an explicit $t$-nomial vanishing 
on $q^{\frac{t-2}{t-1}}$ distinct cosets of $\F^*_q$.   
Over prime fields $\F_p$, computational data we provide suggests that it is 
harder to construct explicit sparse polynomials with many roots. 
Nevertheless, assuming the Generalized Riemann Hypothesis, 
we find explicit trinomials having 
$\Omega\!\left(\frac{\log p}{\log \log p}\right)$ distinct roots in $\F_p$.  
\end{abstract} 

\maketitle 

\vspace{-.5cm}
\section{Introduction}  
How can one best bound the complexity of an algebraic set in terms of the 
complexity of its defining polynomials? Over the complex numbers (or any 
algebraically closed field),
B\'ezout's Theorem \cite{bezout} bounds the number of roots, for a system of 
multivariate 
polynomials, in terms of the degrees of the polynomials. Over finite fields, 
Weil's famous mid-20$\thth$ century result \cite{weil} bounds the number of 
points on a curve in terms of the genus of the curve (which can also 
be bounded in terms of degree). These bounds are 
optimal for dense polynomials. For sparse polynomials, over fields that are 
not algebraically closed, these bounds can be much larger than necessary. For 
example, Descartes' Rule \cite{descartes} tells us that a univariate real 
polynomial with exactly $t$ monomial terms always has 
less than $2t$ real roots, even though the terms may have arbitrarily large 
degree.

Is there an analogue of Descartes' Rule over finite fields? 
Despite the wealth of beautiful and deep 20$\thth$-century results on 
point-counting for curves and higher-dimensional varieties over finite fields, 
varieties defined by sparse {\em univariate} polynomials were all but ignored 
until \cite{cfklls} (see Lemma 7 there, in particular). Aside from their own 
intrinsic interest, refined univariate root counts over finite fields are 
useful in applications 
such as cryptography (see, e.g., \cite{cfklls}), the efficient generation 
of pseudo-random sequences (see, e.g., \cite{bomb}), and refined 
estimates for certain exponential sums over finite fields \cite[Proof of 
Theorem 4]{bourgain}. 
For instance, estimates on the number of roots of univariate tetranomials over 
a finite field were a key step in establishing the uniformity of the 
{\em Diffie-Helman 
distribution} \cite[Proof of Thm.\ 8, Sec.\ 4]{cfklls} --- a quantitative 
statement important in justifying the security of cryptosystems based on the 
Discrete Logarithm Problem.  
 
We are thus interested in the number of roots of sparse univariate 
polynomials over 
finite fields. The polynomial  $x^q-x$ having two terms and exactly $q$ 
roots in $\F_q$ might suggest that there is no finite field analogue of 
Descartes' rule. However, the roots of $x^q-x$ consist of $0$ and the roots of 
$x^{q-1}-1$, and the latter roots form the unit group $\F^*_q\!:=\!\F\setminus
\{0\}$. For an arbitrary binomial $ax^n + bx^m \in \F_q[x]$ with $n\!>\!m$ 
and $a$ and $b$ nonzero, the roots consist of $0$ (if $m>0$) and the roots of 
$x^{n-m} + b/a$. Note that the number of roots of $x^{n-m} + b/a$ in $\F_q$ is 
either $0$ or $\gcd(n-m,q-1)$. In the latter case, the roots form a coset of 
a subgroup of $\F^*_q$. For polynomials with three 
or more terms, the number of roots quickly becomes mysterious and difficult,
 and, as we shall demonstrate in this paper, may exhibit very different 
behaviors in the two extreme cases where (a) $q$ is a large power of a prime,  
and (b) $q$ is a large prime.
 
To fix notation, we call a polynomial 
$f(x)\!=\!c_1+c_2x^{e_2}+\cdots+c_tx^{e_t}\!\in\!\F_q[x]$ with\linebreak  
$e_1\!<\!e_2\!<\cdots<\!e_t\!<\!q-1$ and $c_i\neq 0$ for all $i$  
a {\em (univariate) $t$-nomial}. The best current upper bounds 
on the number of roots of $f$ in $\F_q$, as a function of $q$, $t$, and 
the coset structure of the roots of $f$, can 
be summarized as follows, using $|\cdot|$ for set cardinality:  
\begin{thm} \label{thm:z} 
Let $f\!\in\!\F_q[x]$ be any 
univariate $t$-nomial with degree $<\!q-1$ and 
exponent set $\{e_1,\ldots,e_t\}$ containing $0$. Set 
$\delta(f)\!:=\!\gcd(e_1,\ldots,e_t,q-1)$, 
$Z(f)\!:=\!\left\{x\!\in\!\F_q\; | \; f(x)\!=\!0\right\}$, 
$R(f)\!:=\!\left|Z(f)\right|$, and let $C(f)$ denote the maximum cardinality 
of any coset (of any subgroup of $\F^*_q)$ contained in $Z(f)$. Then: \\ 
0.\ (Special case of \cite[Thm.\ 1]{karpshpar}) 
$R(f)\!\leq\!\frac{t-1}{t}(q-1)$.\\  
1.\ \cite[Thm.\ 1.1]{bcr2} $Z(f)$ is a union of 
no more than $2\left(\frac{q-1}{\delta(f)} \right)^{\frac{t-2}{t-1}}$ cosets, 
each associated\linebreak 
\mbox{}\hspace{.5cm}to one of two subgroups $H_1\!\subseteq\!H_2$ of $\F^*_q$, 
where $|H_1|\!=\!\delta(f)$,  
$|H_2|\!\geq\!\delta(f)\left(\frac{q-1}{\delta(f)}\right)^{1/(t-1)}$, 
and\linebreak 
\mbox{}\hspace{.5cm}$|H_2|$ can be determined within $2^{O(t)} 
(\log q)^{O(1)}$ bit operations.\\  
2. \cite[Thm.\ 1.2]{ko} For $t\!=\!3$ we have 
$R(f)\!\leq\!\delta(f)\left\lfloor 
\frac{1}{2}+\sqrt{\frac{q-1}{\delta(f)}}\right\rfloor$ and, if we have in 
addition\linebreak 
\mbox{}\hspace{.5cm}that $q$ is a square and $\delta(f)\!=\!1$, then 
$R(f)\!\leq\!\sqrt{q}$.\\  
3. (See \cite[Thms.\ 2.2 \& 2.3]{kelley}) For any $t\!\geq\!2$ we have 
$R(f)\!\leq\!2(q-1)^{\frac{t-2}{t-1}}C(f)^{1/(t-1)}$.\\   
\mbox{} \ \ \ \scalebox{.93}[1]{Furthermore, 
$\displaystyle{C(f)\leq\max\{k\!\in\!\N\; : \; 
k|(q-1) \text{ and, for all } i, \text{ there is a } j\!\neq\!i \text{ with } 
k|(e_i-e_j)\}}$. \qed} 
\end{thm}

\noindent 
For any fixed $t\!\geq\!2$, {\em Dirichlet's Theorem} (see, e.g., 
\cite[Thm.\ 8.4.1, Pg.\ 215]{bs}) implies that there are infinitely many prime 
$q$ with $t|(q-1)$. For such pairs $(q,t)$ the bound from Assertion (0) is 
tight: The roots of \[ f(x)\!=\!\frac{x^{q-1}-1}{x^{(q-1)/t}-1}\!=\!1 
+x^{\frac{1}{t}(q-1)}+\cdots+x^{\frac{t-1}{t}(q-1)},\]
are the disjoint union of $t-1$ cosets of size $\delta(f)\!=\!\frac{q-1}{t}$. 
(There are no $H_2$-cosets for this $t$-nomial.) However, Assertions (1) 
and (3) tell us that we can get 
even sharper bounds by making use of the structure of the cosets inside $Z(f)$. 
For instance, when $t\!=\!3$ and $\delta(f)\!=\!1$,  
Assertion (2) yields the upper bound $\sqrt{q}$, which is smaller than 
$\frac{2}{3}(q-1)$ for $q\!\geq\!5$.  
 
While Assertion (3) might sometimes not improve on the upper bound 
$\frac{t-1}{t}(q-1)$, it is often the case that $C(f)$ is provably small enough 
for a significant improvement. For instance, when $e_1\!=\!0$ and 
$\gcd(e_i,q-1)\!=\!1$ for all $i\!\geq\!2$, we have $C(f)\!=\!1$ and 
then $R(f)\!\leq\!2(q-1)^{\frac{t-2}{t-1}}$. 

Our first main result is two explicit families of $t$-nomials revealing that 
Assertions (1)--(3) are close to optimal for {\em non}-prime $q$.  
\begin{thm}  
\label{thm:opt}
Let $t,u,p\!\in\!\N$ with $t\!\geq\!2$ and $p$ prime. 
If $q\!=\!p^{(t-1)u}$ then the polynomial  
\[ r_{t,u,p}(x):=1+x+x^{p^u}+\cdots+x^{p^{(t-2)u}}\]
has $\delta(r_{t,u,p})\!=\!C(r_{t,u,p})\!=\!1$ and exactly 
$q^{(t-2)/(t-1)}$ roots in $\F_q$. Furthermore, if $q\!=\!p^{tu}$, then the 
polynomial 
\[g_{t,u,p}(x):=1+x+x^{1+p^u}+\cdots +x^{1+p^u+\cdots 
+p^{(t-2)u}}\] 
has $\delta(g_{t,u,p})\!=\!1$, $C(g_{t,u,p})\!\leq\!\lfloor t/2\rfloor$, 
and exactly $q^{(t-2)/t}+\cdots+q^{1/t}+1$ roots in $\F_q$.  
\end{thm} 

\noindent 
Theorem \ref{thm:opt} is proved in Section \ref{sec:trick} below.  
The polynomials $r_{t,u,p}$ show that the bounds from Assertions (1) and (3) 
of Theorem \ref{thm:z} {\em are within a factor of $2$ of being optimal}, 
at least for $\delta(f)\!=\!C(f)\!=\!1$ and a positive density of prime powers 
$q$. Note in particular that $r_{3,u,p}$ shows that Assertion (2) of 
Theorem \ref{thm:z} is {\em optimal} for square $q$ and $\delta(f)\!=\!1$.  
(See \cite[Thm.\ 1.3]{ko} for a different set of extremal trinomials 
when $q$ is an odd square.)  
The second family $g_{t,u,p}$ reveals similar behavior for a different family 
of prime powers.

\medskip
Optimally bounding the maximal number of roots (or cosets of roots)  
for the case of {\em prime} $q$ is already more subtle in the trinomial 
case: We are unaware of any earlier lower bound that is a strictly increasing 
function of $q$. (Note, for instance, that $1+x^{(q-1)/3}+x^{2(q-1)/3}$ 
vanishes on just \underline{$\pmb{2}$} cosets of roots, assuming 
$q\!=\!1$ mod $3$.) 
In Section \ref{sec:harder}, we shall prove the following 
theorem.
\begin{thm} \label{thm:grh} 
For any $n\!\geq\!2$ and prime $p\!\geq\!n+2$, consider 
$h_{n,p}(x)\!:=\!x^n-x-1$ as an element of $\F_p[x]$. Then  
$\delta(h_{n,p})\!=\!C(h_{n,p})\!=\!1$, and 
there exists an infinite family of primes $p$ such that the trinomial $h_{n,p}$ 
has exactly $n$ distinct roots in $\F_p$ where 
\begin{enumerate}
\item[(a)] $n\!=\!\Omega\!\left(\frac{\log\log p}{\log\log \log p} \right)$ 
unconditionally, and  
\item[(b)] $n\!=\!\Omega\!\left(\frac{\log p}{\log \log p}\right)$ 
if the Generalized Riemann Hypothesis (GRH) is true. 
\end{enumerate}
\end{thm} 

\noindent 
This theorem is proved by using results from algebraic number theory on the splitting of primes in number 
fields. In particular, we use a classic estimate of Lagarias, Montgomery, and 
Odlyzko \cite{lmo} (reviewed in Section \ref{sec:harder} below) on the 
least prime ideal possessing a Frobenius element exhibiting a particular Galois 
conjugacy class. The latter result is in turn heavily based on the 
effective Chebotarev Density Theorem of Lagarias and Odlyzko \cite{lo}. 

Cohen \cite{Coh70,Coh72} applied an effective function field analogue of 
the Chebotarev Density Theorem to trinomials of the form $x^n + x + b$, where 
$b$ varies in $\F_p$ and $p$ goes to infinity. His results tell us that
when $p$ is large, there always exist $b\!\in\!\F_p$ so that  $x^n + x + b$ 
splits completely over $\F_p$, and the least such $p$ satisfies 
$n\!=\!\Omega\!\left(\frac{\log p}{\log \log p}\right)$ unconditionally. This 
gives the existence of $b\in \F_p$ (and infinitely many primes $p$) 
with $x^n+x+b$ having many roots in $\F_p$, but gives no information on 
how to easily find such $b$. The novelty of Theorem \ref{thm:grh}  
is thus an even simpler family of trinomials with number of roots in  
$\F_p$ increasing as $p$ runs over an infinite sequence of primes. 
Assuming GRH, the rate of increase matches Cohen's lower bound.  

Could the existence of trinomials $f$ with $\delta(f)\!=\!1$ and, 
say, $\Omega\!\left(\sqrt{p}\right)$ roots in $\F_p$ (as one might conjecture 
from Theorem \ref{thm:opt}) be obstructed by $p$ not being a square? 
We feel that $p$ being prime is such an obstruction and, 
based on experimental results below, we conjecture the following upper bound: 
\begin{conj}
\label{conj:logp} 
There is an absolute constant $c\!\geq\!1$ such that, for any prime  
$p\!\geq\!3$ and $\gamma,e_2,e_3\!\in\!\{1,\ldots,p-2\}$, 
with $e_3\!>\!e_2\!>\!0$ and $\gcd(e_2,e_3,p-1)\!=\!1$, 
the trinomial $\gamma+x^{e_2}+x^{e_3}$ has no more than $(\log p)^c$ roots 
in $\F_p$.
\end{conj}

\noindent 
(See also \cite[Conj.\ 1.5 \& Sec.\ 4]{ko} for other refined conjectures 
and heuristics in this direction.) 
It is a simple exercise to show that, to compute the maximal number of roots 
in $\F_p$ of trinomials with $\delta(f)\!=\!1$, one can restrict to the 
family of trinomials in the conjecture.   

For any $n\!\in\!\N$, let $p_n$ be the least prime for which 
there exists a univariate trinomial $f_n$ with 
$\delta(f_n)\!=\!1$ and exactly 
$n$ distinct roots in $\F_{p_n}$.  We did a computer search to find the values 
of $p_n$ for $1\!\leq\!n\!\leq 16$. They are... \newpage 
\begin{table}[ht]
\scalebox{.95}[1]{
\begin{tabular}{r|cccccccccccccccc}
$n$ & 1 & 2 & 3 & 4 & 5 & 6 & 7 & 8 & 9 & 10 & 11 & 12 & 13 & 14 & 15 & 16 \\ 
 \hline 
$p_n$ & 3 & 5 & 11 & 23 & 47 & 151 & 173 & 349 & 619 & 1201 & 2753 
& 4801 
& 10867 
& 16633 & 71237 & 8581 
\end{tabular}}
\end{table} 

\noindent 
For example, $p_{16}\!=\!8581$ because 
$-364+363x+x^{2729}$ has exactly $16$ roots in $\F_{8581}$, and 
{\em no} trinomial $f\!\in\!\F_p[x]$ with 
$p\!<\!8581$ and $\delta(f)\!=\!1$ has more than $15$ roots in $\F_p$. In the 
appendix, we give representative trinomials for each $p_n$ above.  

To get a feel for how the maximal number of roots of a trinomial grows 
with the field size, let us compare the graphs 
(drawn darker) of the functions $0.91\log x$ and $1.77\log x$
with the\linebreak 
\scalebox{.95}[1]{piecewise linear curve (drawn lighter) going through the 
sequence of points $((p_1,1),\ldots,(p_{12},12)$,}\linebreak 
$(p_{16},16),(p_{13},13), (p_{14},14), (p_{15},15))$ as shown in Figure 
\ref{Fig1} below.   
\begin{figure}[ht]
\epsfig{file=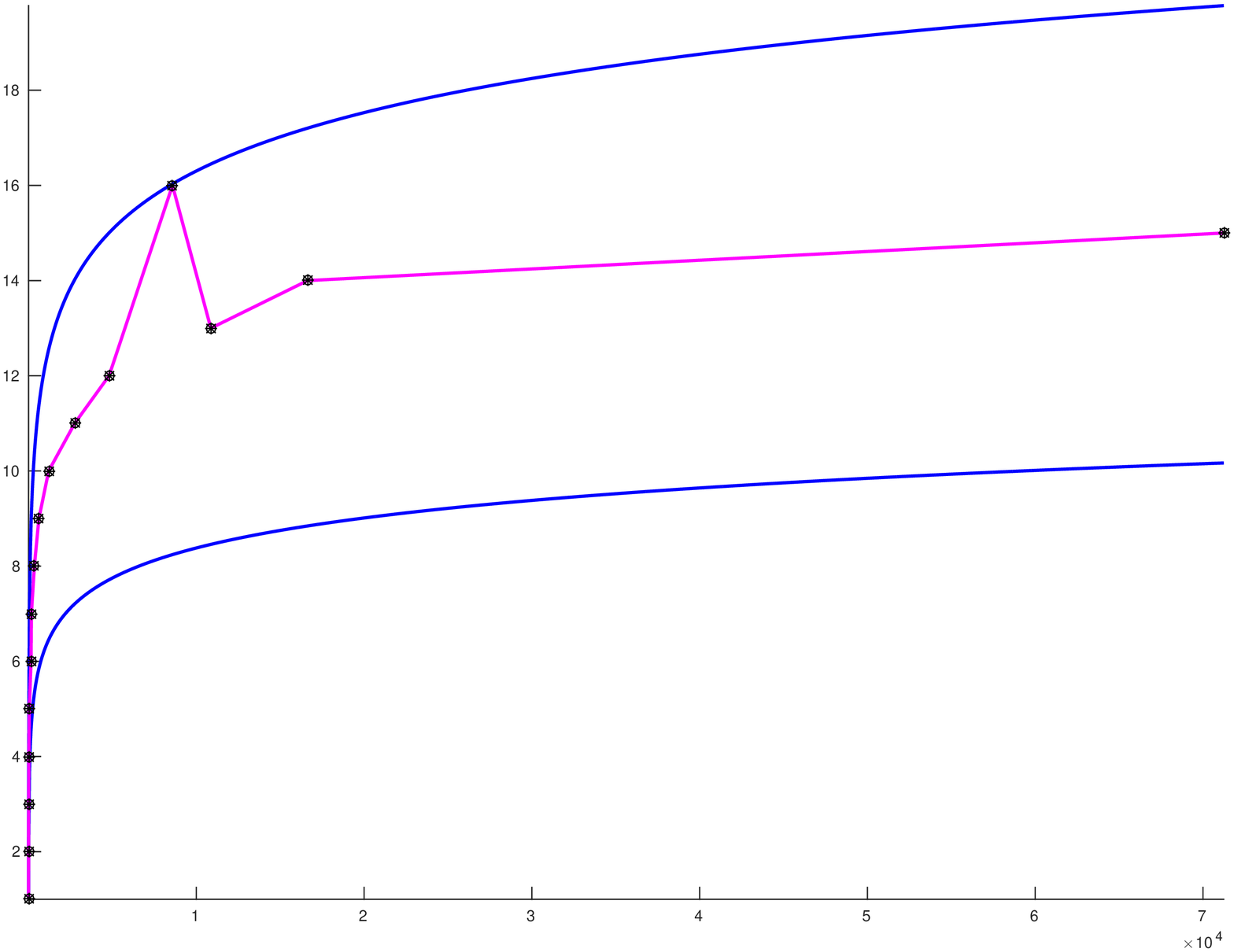,height=3in,clip=} 
\caption{}
\label{Fig1}
\end{figure} 
We used some simple {\tt Maple} and 
{\tt Sage} programs that took a few hours to find most of the data above. 
The value of $p_{15}$ took {\tt C} code (written by Zander 
Kelley) running on a supercomputer for 2 days. 
 
Quantitative results on sparse polynomials over finite fields sometimes 
admit stronger analogues involving complex roots of unity. For instance, 
\cite{bcr1,bcr2} and \cite{cheng} deal with the complexity of deciding whether 
a sparse polynomial vanishes on a (finite) subgroup of, respectively, 
$\F^*_q$ or $\C^*$. It is thus interesting to observe a recent complex analogue 
to our trinomial root counts over $\F_q$: 
Theobald and de Wolff \cite{tt} proved that, if $\gcd(e_2,e_3)\!=\!1$, a 
trinomial $c_1+c_2x^{e_2}+c_3x^{e_3}\!\in\!\C[x]$ can have at most  
$2$ complex roots of the same absolute value. So any 
such trinomial has at most $2$ roots in the {\em torsion subgroup} 
$\{\zeta\!\in\!\C\; | \; \zeta^n\!=\!1 \text{ for some } 
n\!\in\!\N\}$ of $\C^*$. This upper bound is sharp: Consider 
$(x-1)(x-\zeta)$ for any $\zeta\!\not\in\!\{\pm 1\}$ satisfying 
$\zeta^n\!=\!1$ for some $n\!\geq\!3$.

\section{Main Lower Bounds in the Prime Power Case} 
\label{sec:trick}  

\noindent 
{\bf Proof of Theorem \ref{thm:opt}}: 
To establish the root count for $r_{t,u,p}$ it clearly suffices 
to prove that $r_{t,u,p}$ divides $x^q-x$ or, equivalently, 
$x^{p^{(t-1)u}}\!=\!x$ in the ring 
$R\!:=\!\F_p[x]/\langle r_{t,u,p}(x)\rangle$.  
\linebreak  
Toward this end, observe that $x^{p^{(t-1)u}} = 
\left(x^{p^{(t-2)u}} \right)^{p^u} = 
\left(-1-x^{p^0}-\cdots-x^{p^{(t-3)u}}\right)^{p^u}$ in $R$. Since 
$(a+b)^{p^u}\!=\!a^{p^u}+b^{p^u}$ in any ring of characteristic 
$p$, we thus obtain by induction on $t$ that 
$x^{p^{(t-1)u}}=(-1)^{p^u}\left(1+x^{p^u}+\cdots+p^{(t-2)u}\right)$ 
in $R$. The last factor is merely $r_{t,u,p}(x)-x$, so we obtain 
$x^{p^{(t-1)u}}=(-1)^{p^u}(-x) = (-1)^{1+p^u}x$ 
in $R$. Since $(-1)^{1+p^k}\!=\!1$ in $\F_p$ for all primes $p$, 
we have thus established the root count for $r_{t,u,p}$. 

That $\delta(r_{t,u,p})\!=\!1$ is clear since 
$r_{t,u,p}$ has a nonzero constant term and $1$ as one of its exponents. 
Likewise, $\delta(g_{t,u,p})\!=\!1$. 
That $C(r_{t,u,p})\!=\!1$ is clear because the lowest 
exponents of $r_{t,u,p}$ are $0$ and $1$, the rest are powers of $p$, 
and $\gcd(p^k,q-1)\!=\!1$ for all $k\!\in\!\N$. 
We postpone proving our upper bound on $C(g_{t,u,p})$ until after we prove 
our stated root count for $g_{t,u,p}$. 

Consider now the set $S$ of elements in $\F_{q}$ whose trace to $\F_{p^u}$ 
is zero, that is,
\[ S:=\left\{ a \in \F_q \; | \;  a + a^{p^u} + a^{p^{2u}} + \cdots +  
a^{p^{(t-1)u}} = 0\right\}.\]
Then $S$ has  $q/p^u\!=\!p^{(t-1)u}$ elements and is a vector space of 
dimension $t-1$ over $\F_{p^u}$. Let $a \in S$ be nonzero. We show that 
$a^{p^u -1}$ is a root of $g_{t,u,p}$: 
\[ g_{t,u,p}(a^{p^u -1})= 1 +  \sum_{i=1}^{t-2} a^{(p^u -1) (p^{iu} + \cdots 
+ p^u + 1)} = 1 + \sum_{i=1}^{t-2} a^{p^{(i+1)u} - 1} = 
\frac{1}{a} \sum_{i=0}^{t-1}  a^{p^{iu}} =0.\]
Now note that, for any $a\!\in\!S$ and any nonzero $w\!\in\!\F_{p^u}$, the 
element $aw$ is also in $S$. Also, for $a,b\!\in\!\F_q$, 
$a^{p^u-1}\!=\!b^{p^u-1}$ if and only if $b\!=\!aw$ for some nonzero 
$w\!\in\!\F_{p^u}$. Therefore, when $a$ runs through  $S\!\setminus\!\{0\}$, 
the element $a^{p^u-1}$ yields $(p^{(t-1)u} - 1)/(p^u -1)\!=\! 
1 + p^u + \cdots + p^{(t-2)u}$ roots for $g_{t,u,p}$. 

To finally prove our upper bound on $C(g_{t,u,p})$, note that, for $j\!>\!i$, 
\[1+p^u+\cdots+p^{(j-2)u}\!=\!\frac{p^{(j-1)u}-1}{p^u-1},\]
and
\[ \frac{p^{(j-1)u}-1}{p^u-1}-\frac{p^{(i-1)u}-1}{p^u-1} 
\!=\!p^{(i-1)u}\left(\frac{p^{(j-i)u}-1}{p^u-1}\right).\]
So for $i\!\geq\!2$,
\begin{eqnarray*}
\max\limits_{j\in\{1,\ldots,t\}\setminus\{i\}}  \gcd\!\left(p^{(i-1)u}\left(\frac{p^{(j-i)u}-1}{p^u-1}\right), p^{tu}-1\right)
 & = & \max\limits_{j\in\{1,\ldots,t\}\setminus\{i\}}  \gcd\!\left(\frac{p^{|j-i|u}-1}{p^u-1},p^{tu}-1\right)\\
  & = & \max\limits_{\ell\in\{1,\ldots,\max\{i-1,t-i\}\}}  \gcd\!\left(\frac{p^{\ell u}-1}{p^u-1},p^{tu}-1\right).
\end{eqnarray*}
Hence, \ \ $D(g_{t,u,p}) := 
\min\limits_{i\in\{1,\ldots,t\}} \max\limits_{\ell\in\{1,\ldots,\max\{i-1,t-i\}\}}  \gcd\!\left(\frac{p^{\ell u}-1}{p^u-1},p^{tu}-1\right)$
\begin{eqnarray*}
& = & \max\limits_{\ell\in\{1,\ldots,\lfloor t/2\rfloor \}\}} \gcd\!\left(\frac{p^{\ell u}-1}{p^u-1},p^{tu}-1\right)\\
& = &  \max\limits_{\ell\in\{1,\ldots,\lfloor t/2\rfloor \}\}} \gcd\!\left(\frac{p^{\ell u}-1}{p^u-1},\left(\frac{p^{tu}-1}{p^u-1}\right) (p^u-1)\right)\\ 
& \leq &   \max\limits_{\ell\in\{1,\ldots,\lfloor t/2\rfloor \}\}} \gcd\!\left(\frac{p^{\ell u}-1}{p^u-1},\frac{p^{tu}-1}{p^u-1}\right)
                \gcd\!\left(\frac{p^{\ell u}-1}{p^u-1},p^u-1\right),\\ 
 & = & \max\limits_{\ell\in\{1,\ldots,\lfloor t/2\rfloor \}\}} \left(\frac{p^{\gcd(\ell,t) u}-1}{p^u-1}\right)\gcd\!\left(\ell,p^u-1\right).
 \end{eqnarray*}
The last equality follows easily from two elementary facts: 
(1) $\gcd(x^\ell-1,x^t-1)\!=\!x^{\gcd(\ell,t)}-1$,\linebreak 
and (2) $(x-1)(x^{\ell-2}+2x^{\ell-3}+\cdots+(\ell-2)x
+(\ell-1))\!=\!  x^{\ell-1}+\cdots+x^2+x-(\ell-1)$.\linebreak  
So $D(g_{t,u,p})\!\leq\!\max\limits_{\ell\in\{1,\ldots,\lfloor t/2\rfloor \}} 
1\cdot \ell \leq \lfloor t/2\rfloor$. By \cite[Prop.\ 2.4]{kelley} and 
\cite[Thm.\ 2.2]{kelley},\linebreak $D(g_{t,u,p})\!\geq\!C(g_{t,u,p})$, so we 
are done. \qed

\section{Main Lower Bounds in the Prime Case}  \label{sec:harder} 
We'll need several results from algebraic number theory. First, let $K$ be any 
number field, i.e., a finite algebraic extension of $\Q$. Let  $d_{K}$ denotes 
the discriminant of $K$ over $\Q$, and $\cO_K$ the ring of algebraic integers 
of $K$, i.e., 
those elements of $K$ with monic minimal polynomial in $\Z[x]$.  We need to know the size of the smallest 
prime $p \in \Z$ that splits completely in $\cO_K$.    There are various bounds in the literature that are proved via some effective version of the 
Chebotarev Density Theorem. For instance:  
\begin{thm} 
\label{thm:least} 
(See \cite[Cor.\ 1.2 \& pp.\ 461--462]{lo} and \cite[Thm.\ 1.1]{lmo}.)   
If $f\!\in\!\Z[x]$ is any irreducible polynomial of degree $n$ then 
the least prime $p$ for which the reduction of $f$ mod $p$ has 
$n$ distinct roots in $\F_p$ is (unconditionally) no greater than 
$d^{O(1)}_K$, where $K\!\subset\!\C$ is the splitting field of $f$. 
Furthermore, if GRH is true, then $p\!=\!O((\log d_K)^2)$. 
\end{thm} 

\noindent 
The papers \cite{lo,lmo} in fact work in much greater generality: Our 
limited paraphrase above is simply the special case where 
one is looking for a prime yielding a Frobenius element corresponding 
to the identity element of the Galois group of $f$ over $\Q$. 

The best recent estimates indicate that, in the 
unconditional case of Theorem \ref{thm:least}, we can take the $O$-constant to 
be $40$, for sufficiently large $d_K$ \cite{kadiring}.  
Also, for an abelian extension $K$ over $\Q$, Pollack \cite{Poll} gives a 
much better bound (in the unconditional case): 
$p\!=\!O_{\eps,K}\!\left(d^{\frac{1}{4}+ \eps}_K\right)$ where 
$\eps>0$ is arbitrary, and the implied $O$-constant depends only on $\eps$ and 
the degree of $K$ over $\Q$. 

We will also need good bounds on discriminants of number fields. In the 
following 
theorem, the lower bound is  due to Minkowski \cite{minkowski} and the 
upper bound  follows easily from work of T\^oyama \cite{toyama} (by induction 
on the number of composita generated by the distinct roots of $f$). 
\begin{thm} 
(See, e.g., \cite[pp.\ 259--260]{bs}.) 
\label{thm:disc} 
For any number field $K$ of degree $n$ over $\Q$,  we have $
d_K\!\geq\!\frac{n^{2n}}{(n!)^2}\!\geq\! 
\frac{(\pi e^2/4)^n}{2\pi n}\!>\!\frac{5.8^n}{6.3n}$.  
Also,  if $K$ has minimal polynomial  $f\!\in\!\Q[x]$ and $L$ is the splitting 
field of $f$, then $d_L$ divides 
$d_K^{(n-1)!+(n-2)! n +\cdots+1! n^{n-2}+0!n^{n-1}}$. \qed 
\end{thm} 

\noindent 
{\bf Proof of Theorem \ref{thm:grh}}: Clearly, $\delta(h_{n,p})\!=\!1$. 
Also, since $\gcd(n,n-1)\!=\!1$, it is clear that $C(h_{n,p})\!=\!1$. 
Now, for any $n \geq 2$, the trinomial 
$h_n:=x^n -x -1$ is irreducible over $\Q$ \cite{selmer}.
Let $\alpha \in \C$ be any root of $h_n$ and let $K = \Q(\alpha)$, so that 
$[K:\Q]=n$. Then $d_K$ divides the resultant of $h_n$ and $h'_n$ 
\cite[Thm.\ 8.7.1, pg.\ 228]{bs}. The resultant of 
$h_n$ and $h'_n$ can then be computed explicitly to be 
$(-1)^{\frac{(n+2)(n-1)}{2}} \left(n^n+(-1)^n (n-1)^{n-1}\right)$ 
\cite{swan}. 
Hence $d_K$ divides $n^n+(-1)^n (n-1)^{n-1}$.
(We note that an elegant modern development of trinomial discriminants 
can be found in Chapter 12 of \cite{gkz94}; see Formula 
1.38 on Page 406 in particular).

Let $L$ be the the splitting field of $h_n$. Then $L$ has degree at most $n!$ 
and, by Theorem \ref{thm:disc}, 
$d_K\!>\!\frac{5.8^n}{6.3n}$ and $d_L$ divides $\left(n^n+ (-1)^n (n-1)^{n-1}\right)^{(n-1)!+(n-2)! n + \cdots+1! n^{n-2}+0! n^{n-1}}$. 
Note that, for $n \geq 3$, we have  
$n^n+(-1)^n (n-1)^{n-1}\! \leq n^n+ (n-1)^{n-1}\!  \leq\!e^{n\log n + \frac{4}{27}}$.

Also, by Stirling's Estimate \cite[Pg.\ 200]{rudin}, 
$n!\!<\!e\sqrt{n}\left(\frac{n}{e}\right)^n$ (for all  $n\!\geq\!1$), so we have
\begin{eqnarray*} 
& &  (n-1)!    +n(n-2)!+\cdots+2!n^{n-3}+1!n^{n-2}+0!n^{n-1}\\ 
& &\ \  <  \ e\sqrt{n-1}\left(\frac{n-1}{e}\right)^{n-1}+   e\sqrt{n-2}\left(\frac{n-2}{e}\right)^{n-2}n+  
\cdots+e\sqrt{1}\left(\frac{1}{e}\right)^1 n^{n-2}+ n^{n-1}\\ 
& & \ \ < \   e\sqrt{n}\left(1+\frac{1}{e}+\cdots+\frac{1}{e^{n-1}}\right)n^{n-1}
\!=\!\left(\frac{e}{1-1/e}\right)n^{n-1/2}\!<\!4.31n^{n-1/2},
\end{eqnarray*}
and thus $d_L\!<\!e^{4.5n^{n+1/2}\log n}$.  

Theorem \ref{thm:least} then tells us that there is a prime $p \in \Z$ so that 
$h_n$ splits completely modulo $p$ with no repeated roots where
\begin{enumerate}
\item[(a)] $p=e^{O(n^{n+1/2}\log n)}=  e^{e^{(n+1/2+o(1))\log n}}$ 
unconditionally, and  
\item[(b)] $p=O((n^{n+1/2}\log n)^2)=e^{(2n+1+o(1))\log n}$ if 
GRH is true.
\end{enumerate}
So then we must have
\begin{enumerate}
\item[(a)]  $(n+1/2+o(1))\log n\!\geq\!\log\log p$ unconditionally, and 
\item[(b)] $(2n+1+o(1))\log n\!\geq\!\log p$ if GRH is true. 
\end{enumerate}
Considering the growth order of the inverse function of 
$x\log x$, we obtain our stated asymptotic lower bounds in Theorem 
\ref{thm:grh}. \qed 

\medskip
We used the family of trinomials $x^n-x-1$ mainly for the sake of 
simplifying our proof. Many other families would likely exhibit the same 
behavior, albeit with some additional intricacies in applying prime ideal 
bounds. However, the deeper question is to understand the structure of 
{\em truly} extremal trinomials over prime fields, such as those appearing 
in the Appendix below. 

\section*{Acknowledgements}
We thank the American Institute of Mathematics for their
splendid hospitality and support of our AIM SQuaRE project, which
formed the genesis of this paper. Special thanks go to Zander Kelley
(and the Texas A\&{}M Supercomputer facility) for computing $p_{15}$,
and pointing out the paper \cite{karpshpar}.
Thanks also to Kiran Kedlaya for sharing his cython code (which helped
push our computational experiments further) and for pointing out an error
in an earlier computation of $p_{12}$.
We also thank E.\ Mehmet Kiral, Timo de Wolff, and Matthew P.\ Young for
useful discussions, and the referees for their insightful comments.

\section*{Appendix: Some Extremal Trinomials} 
We list in Figure \ref{Fig2}, for $n\!\in\!\{1,\ldots,16\}$, trinomials $f_n$ 
with $\delta(f_n)\!=\!1$ and $f_n$ having exactly $n$ distinct roots in 
$\F_{p_n}$, with $p_n$ the smallest prime admitting such a trinomial.  
\begin{figure}[ht]
\begin{tabular}{r|c|l} 
$n$ & $f_n$ & $p_n$ \\ \hline 
$1$ & $1+x-2x^2$ & $3$\\
$2$ & $1+x-2x^2$ & $5$\\ 
$3$ & $1-3x+2x^3$ & $11$\\ 
$4$ & $-2+x+x^4$ & $23$\\       
$5$ & $1+4x-5x^8$ & $47$\\ 
$6$ & $1+24x-25x^{33}$ & $151$\\ 
$7$ & $-2+x+x^{34}$ & $173$ \\ 
$8$ & $1+23x-24x^{21}$ & $349$\\ 
$9$ & $-71+70x+x^{184}$ & $619$\\
$10$ & $1+5x-6x^{152}$ & $1201$\\
$11$ & $-797+796x+x^{67}$ & $2753$\\ 
$12$ & $-82+81x+x^{1318}$ & $4801$\\ 
$13$ & $-1226+1225x+x^{225}$ & $10867$\\ 
$14$ & $-39+38x+x^{2264}$ & $16633$ \\  
$15$ & $29574-29573x-x^{27103}$ & $71237$ \\ 
$16$ & $-364+363x+x^{2729}$ & $8581$ 
\end{tabular}
\caption{Trinomials with exactly $n$ distinct roots in $\F_{p_n}$ and $p_n$ 
minimal}
\label{Fig2}
\end{figure}
In particular, for each $n\!\in\!\{1,\ldots,16\}$, a full search was 
done so that the trinomial $f_n$ below    
has the least degree among all trinomials over $\F_{p_n}$ having 
exactly $n$ roots in $\F_{p_n}$. (It happens to be the case that, 
for $n\!\in\!\{1,\ldots,16\}$, we can also pick the middle degree 
monomial of $f_n$ to be $x$.) By rescaling the variable as necessary, we have 
forced $1$ to be among the roots of each of the trinomials below.     
It is easily checked via the last part of Assertion (3) of Theorem 
\ref{thm:z} that $C(f_n)\!=\!1$ for each $n\!\in\!\{1,\ldots,16\}$.

The least prime $p_{17}$ for which there is a trinomial $f_{17}$ with 
$\delta(f_{17})\!=\!1$ and exactly $17$ roots in $\F_{p_{17}}$ is currently 
unknown (as of July 2016). Better and faster code should hopefully change this 
situation soon. 

\bibliographystyle{amsalpha}

\end{document}